\documentclass[11pt]{amsart}
\usepackage{amssymb}
\usepackage{latexsym}
\usepackage{graphicx}

\newtheorem{theorem}{Theorem}
\newtheorem{lemma}{Lemma}

\theoremstyle{definition}

\theoremstyle{remark}

\numberwithin{equation}{section}

\newcommand {\SC} {{\mathbb C}}
\newcommand {\SD} {{\mathbb D}}
\newcommand {\DD} {{\mathbb D}}

\newcommand {\SR} {{\mathbb R}}

\newcommand{\RR}{{\mathbb R}}
\newcommand{\B}{{\bf B}}

\begin{document}

\title[Brownian motion in domains in Euclidean space]{On the duration of stays of Brownian motion in domains in Euclidean space}


\date{May 2021}

\author[D. Betsakos]{Dimitrios Betsakos}
\address{D. Betsakos: Department of Mathematics, Aristotle University of Thessaloniki, GR-54124 Thessaloniki, Greece.} \email{ betsakos@math.auth.gr}

\author[M. Boudabra]{Maher Boudabra}
\address{M. Boudabra: Department of Mathematics, Monash University, Australia} \email{maher.boudabra@monash.edu}

\author[G. Markowsky]{Greg Markowsky}
\address{G. Markowsky: Department of Mathematics, Monash University, Australia} \email{greg.markowsky@monash.edu}


\begin{abstract}
Let $T_D$ denote the first exit time of a Brownian motion
from a domain $D$ in $\RR^n$. Given domains $U,W \subseteq \RR^n$ containing the origin,
we investigate the cases in which we are more likely to have fast exits from $U$ than $W$, meaning ${\bf P}(T_U<t) > {\bf P}(T_W<t)$ for $t$ small. We show that the primary factor in the probability of fast exits from domains is the proximity of the closest regular part of the boundary to the origin.  We also prove a result on the complementary question of longs stays, meaning ${\bf P}(T_U>t) > {\bf P}(T_W>t)$ for $t$ large. This result, which applies only in two dimensions, shows that the unit disk has the lowest probability of long stays amongst all Schlicht domains.
\end{abstract}

\maketitle


\bigskip

\section{{\bf Introduction}}

The distribution of the exit time of a Brownian motion from
a domain in $\RR^n$ gives a measure of the size and information about the shape of the domain. Naturally, a small domain will have an exit time which is smaller, in some sense, than that of a large domain. However, making this statement precise is a bit tricky, and this question has been of interest to a number of researchers, especially in two dimensions. For examples, the reader is referred to the recent papers \cite{Bancar, Bancar2, Ban, Kim, Hern}, as well as to the older works \cite{burk, D}. This paper is concerned with this question, and particularly with the relationship between the probability that the exit time is small and the proximity of the boundary of the domain to the starting point of the Brownian motion.

\medskip

Before stating our results, let us fix notation. Let $\B_t$, $t\geq 0$ denote standard Brownian motion
moving in $\RR^n$, $n\geq 2$. We denote by ${\bf P}$ and ${\bf E}$ the
corresponding probability measure and expectation, respectively, and will use superscripts in ${\bf P}^x$ and ${\bf E}^x$ to signify that we are using the probability measure associated with Brownian motion starting from $x$ (see e.g. \cite{PS}). If ever the starting point is not mentioned, then we assume it is the origin.
For a domain $D$ containing the origin in $\RR^n$, we denote by $T_D$
the first exit time of $\B_t$ from $D$; that is
$$
T_D=\inf\{t \geq 0: \B_t\notin D\}.
$$

We also let ${\rm d}(A)$ denote the distance from the origin to the set $A\subset \SR^n$; that is,
$$
{\rm d}(A) = \inf\{\|x\|: x \in A\}.
$$

For our first result, we consider the quantity ${\bf P}(T_D < t)$ for $t$ small (this is what we mean by ``fast exits") and its relation to ${\rm d}(\partial D)$. The authors of this paper have already proved some results in this direction, in \cite{BBM}. There, we considered this problem only in two dimensions, and assumed also that the domains in question were simply connected. Our main result in that paper was the following.

\medskip

{\bf Theorem A.} {\cite{BBM}} \label{T2} 
{\it Suppose that $U,W$ are two simply connected domains in $\RR^2$ both containing the origin, and that ${\rm d}(\partial U) < {\rm d}(\partial W)$. Then, for all sufficiently small $t>0$,
\begin{equation}\label{T2e}
{\bf P}^0(T_U<t)> {\bf P}^0(T_W< t).
\end{equation}
In fact,
$$
\lim_{t \to 0^+} \frac{{\bf P}(T_U<t)}{{\bf P}(T_W< t)} = \infty.
$$
}

The proof there depended heavily upon simple connectivity and the topology of the plane, and did not seem at the time to generalize to higher dimensions. However, subsequent study has shown that the result does admit a considerable generalization to $\RR^n$, and furthermore that the requirement of simple connectivity was unnecessary. Dropping the requirement of simple connectivity, however, does significantly complicate the proof, and compels us to introduce the concept of {\it regular} boundary points, which we now discuss.


For a closed set $K\subset \SR^n$, let $\tau_K = \inf \{t>0: B_t \in K \}$. We will call this the {\it hitting time} of the set $K$. A point $x$ is called \textit{regular} for $K$ if ${\bf{P}}^x (\tau _K=0)=1$ and \textit{irregular} otherwise (see \cite{PS} or \cite{szi}). 
Intuitively speaking, a regular point $x$ for $K$ is a point from which the Brownian motion hits $K \backslash \{x\}$ immediately upon starting at $x$. Note that ${\bf{P}}^x (\tau _K=0) \in \{0,1 \}$ by Blumenthal's zero-one law; that is, if $x$ is irregular then ${\bf{P}}^x (\tau _K=0)=0$.  
Let  $K^{\rm r}$ denote the set of regular points of  $K$. This set contains the interior of $K$ . This leaves only points in $\partial K$ as unknowns, and these may or may not be regular. A simple example of a planar domain with an irregular boundary point is $\DD \backslash \{0\}$, and the boundary point $0$ is irregular. On the other hand, it is well known that all boundary points in any simply connected domain in the plane (strictly smaller than $\SC$ itself) are regular (see \cite{Ran}, for instance).

The main result of this paper can now be stated, and is as follows.

\begin{theorem}\label{T1}
Let $U,W$ be domains in $\RR^n$, both containing the origin.
Suppose that ${\rm d}((\partial U)^r)<{\rm d}((\partial W)^r)$. Then
\begin{equation} \label{T1e1}
\lim_{t\to 0+}\frac{{\bf P}^0(T_U<t)}{{\bf P}^0(T_W<t)}=\infty.
\end{equation}	
\end{theorem} 

\medskip

We also consider the complementary problem of ``long stays"; that is, the behavior of the quantity ${\bf P}(T_D < t)$ for $t$ large. For this problem, we will focus only on the two-dimensional case. Again, we have a previous result on this topic in our previous paper \cite{BBM}. To describe this result, we first need a definition. For a domain $D\subset \SR^2$, let
$$
{\rm H}(D)=\sup \{p > 0: {\bf E}[(T_D)^p] < \infty\};
$$
note that ${\rm H}(D)$ is proved in \cite{burk} to be exactly
equal to half of the Hardy number of $D$, a purely analytic
quantity, as defined in \cite{hansen}, and is therefore calculable for a number of common domains.
The domains we consider satisfy a normalization condition: they are {\it Schlicht}. A planar simply connected doamin $D\subsetneqq \SC$ is Schlicht if it contains the origin and $D=f(\SD)$, where $f$ is the Riemann map with $f(0)=0$ and $f^\prime(0)=1$. It is known that 
 ${\rm H}(D)
\geq \frac{1}{4}$ as long as $D \neq \SC$ is Schlicht (\cite{burk}). The following is our prior result.

\medskip

{\bf Theorem B.} {\cite{BBM}} \label{T3prep} 
{\it Suppose that $U,W$ are Schlicht domains and  ${\rm H}(U) > {\rm H}(W)$. Then
\begin{equation}\label{T3e}
\limsup_{t \to \infty} \frac{{\bf P}(T_W > t)}{{\bf P}(T_U > t)} = \infty.
\end{equation}
}

We conjectured there that this proposition remained true with the $\limsup$ replaced by $\lim$, but we were not able to prove it (except when $W$ is a wedge). Subsequent study has revealed another result in this direction. Before stating the result, we provide a bit of motivation.

In \cite{D}, Davis  explored the relation of planar Brownian motion  to classical complex analysis (for more on this topic, see \cite{M00} and the references therein). Among a number of important ideas was Davis' statement that ``the distribution of
$T_D$ is an intuitively appealing measure of the size of
$D$" (notation changed to match that in this paper). Davis then suggested applying this idea to the set of Schlicht domains. He conjectured that if $D$ is a Schlicht domain, then
$$
 {\bf P}(T_D<t)\leq {\bf P}(T_\DD<t),
$$
for all $t>0$. However, McConnell disproved this for sufficiently small $t$ and $D$ an infinite strip, in \cite{MAC}; note that this follows also from our results on fast exits (Theorem A), since ${\rm d}(\partial D) < 1$ for any Schlicht domain $D$ other than $\DD$, as a consequence of Schwarz's lemma. 
We will prove that Davis' conjecture is correct for large $t$, i.e. for long stays. Our result is as follows.

\begin{theorem}\label{T3}
Let $D$ be a Schlicht domain other than $\DD$. Then there exists $t_o>0$ such that for every $t>t_o$,
\begin{equation} \label{T1e1new}
{\bf P}^0(T_D < t)<{\bf P}^0(T_\DD < t).
\end{equation}
\end{theorem} 

To prove Theorems \ref{T1} and \ref{T3}, we will use a number of known deep results. The next section contains the necessary preliminaries, and the subsequent section contains the proofs. A short final section contains some concluding remarks.



\section{Preliminaries}
In this section, we collect the topics and results we will need for the proofs of our theorems.
\subsection{Killed Brownian motion}\label{KBM}
 A Brownian motion running in $\RR^n$ admits the transition density 
 $$
 p(t,x,y) = \frac{1}{(2\pi t)^{n/2}} e^{\frac{-|x-y|^2}{2t}}.
 $$
  This means that the probability that Brownian motion starting at $x$ is in a Borel set $A$ at time $t$ is equal to $\int_A p(t,x,y)dy$. Note that $p$ satisfies the heat equation: $\partial_t p = \frac{1}{2}\Delta p$. Our primary interest will be in the ``stopped transition density" $p_D(t,x,y)$ (see e.g. \cite{PS}), which is taken in relation to a domain $D$ in $\RR^n$ and applies to Browian motion killed upon leaving $D$. In particular, if $K$ is outside $D$, then $\int_K p_D(t,x,y)dy=0$. By the strong Markov property, the formula of $p_D(t,x,y)$ is given by 
$$
p_D(t,x,y)=p(t,x,y)-{\bf E}^x(p(t-T_D, B_{T_D},y)1_{T_D<t})
$$ 
and one can see that $p_D(t,x,y)\leq p(t,x,y)$. Note that for every Borel set $K\subset \SR^n$, the identity 
$$
{\bf P}^x(B_t \in K, t< T_D)=\int_K p_D(t,x,y)\,m_n(dy)
$$
persists; here and below $m_n$ is the $n$-dimensional Lebesgue measure. Furthermore, $p_D$ also satisfies the heat equation in $D$  with zero boundary values.

\medskip

Two basic properties of the transition density $p_D$ are
the domain monotonicity:
if $D\subset \Omega$ then $p_D(t,x,y)\leq p_\Omega(t,x,y)$,
and the semigroup property:
$$
p_D(t+s,x,y)=\int_D p_D(t,x,z)p_D(s,z,y)\;m_n(dz),\;\;\;t,s>0,\;x,y\in D.
$$

\medskip

An additional property of $p_D$ is the principle of ``not feeling the boundary". We will need this principle in its simplest form (see \cite{C}): If $D$ is a convex domain in $\SR^n$, then for every $x,y\in D$, 
$$
\lim_{t\to 0+} \frac{p_D(t,x,y)}{p(t,x,y)}=1.
$$ 

\medskip

When $D$ is a ball, the stopped transition density has a radial monotonicity property. Probably, this property is known to experts, but we haven't found it in the literature, and thus we provide a proof.

\begin{lemma}\label{L0}
	Let $D$ be the ball in $\SR^n$ centered at the origin with radius $R$. Then for every $t>0$, the function $p_D(t,0,x)$ is radially strictly decreasing.
\end{lemma}	
\proof

\begin{figure}
\begin{centering}
\includegraphics[width=10cm,height=10cm,keepaspectratio]{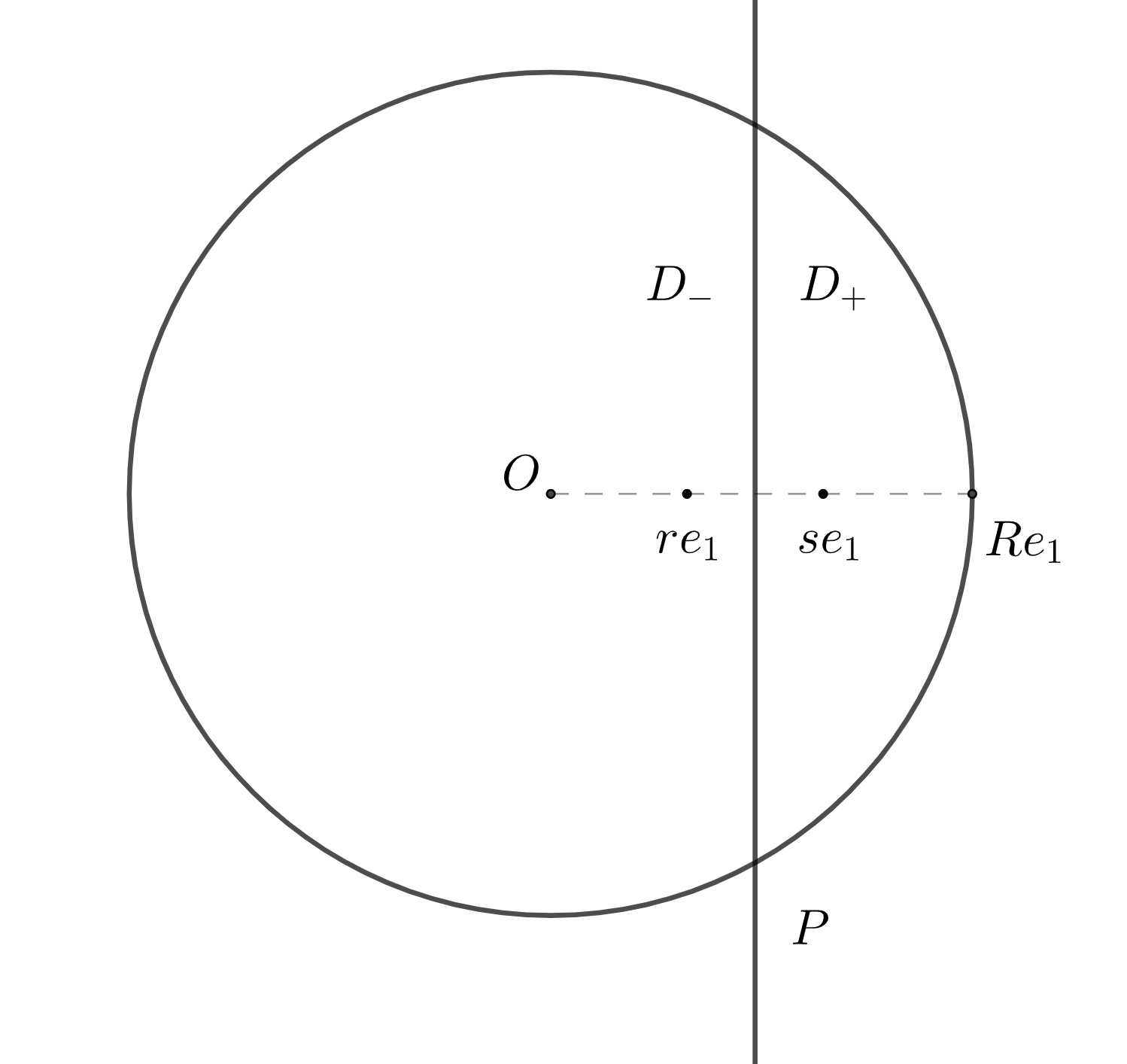}
\par\end{centering}
\caption{The ball, $D_+$, and $D_-$.}
\end{figure} 

By symmetry, $p_D(t,0,x)$ is a radial function. Let $e_1=(1,0,\dots,0)$ be the first coordinate vector
and let $0<r<s<R$. Consider the $(n-1)$-dimensional plane $P$ (a line when $n=2$), perpendicular to $e_1$ and passing from the point $\frac{r+s}{2}\,e_1$. Let 
$$
D_+:=\{x=(x_1,x_2,\dots,x_n)\in D:x_1>\frac{r+s}{2}\},
$$
$$
D_-:=\{x=(x_1,x_2,\dots,x_n)\in D:x_1<\frac{r+s}{2}\}.
$$
Note that the reflection of $D_+$ in $P$ is contained in $D_-$, and that the reflection of the point $se_1$ in $P$ is the point $re_1$. 

Now, since $p_D$ satisfies the heat equation, we can use a special case of a polarization result \cite[Theorem 9.4]{BS}, \cite[Theorem 4]{B2}, and conclude that
\begin{equation}\label{L0p1}
p_D(t,0,se_1)<p_D(t,0,re_1).
\end{equation}
\qed

\medskip

\subsection{Green's function}\label{GF}
An important related object is {\it Green's function}. We give its probabilistic definition; see e.g. \cite{PS}, \cite{szi}. With $p_D$ as above, we define
$$
G_D(x,y) = \int_0^\infty p_D(s,x,y) ds.
$$

\medskip

\subsection{Logarithmic and Newtonian capacity}\label{LNC}
We now discuss the {\it capacity} of a set in $\RR^n$. The definitions of capacity of sets in $\RR^n$ are not entirely standardized in the literature; we will follow the definitions used in \cite{PS}, \cite{Ran}. For a compact set $K$ in $\RR^n$ with $n \geq 2$, define
$$ 
R_n(K) := \inf \int _{K \times K } f(x-y) \mu (dx) \mu (dy),
$$ 
where 
$$
f(x) = \begin{cases}
-\ln|x|, & \text{if }\ensuremath{n=2}\\
|x|^{2-n}, & \text{if \ensuremath{n \geq 3}},
\end{cases}
$$
and where the infinimum is taken over all probability measures $\mu$ having their support on $K$.
If $K\subset \SR^2$, its {\it logarithmic capacity} is defined by $c_2(K)=e^{-R_2(K)}$. If $K\subset \SR^n$, $n\geq 3$, its  {\it Newtonian capacity} is $c_n(K)=R_n(K)^{-1}$. 
These capacities have an important connection with Brownian motion in $n$ dimensions, as we now describe. 

A compact  set $K$ is referred to as \textit{polar} when ${\bf P} (\tau_K < \infty )=0$, i.e  when it is not hit by Brownian motion with probability one. Otherwise, it is called \textit{nonpolar}.
By a theorem of Kakutani (see e.g. \cite[Prop. 6.1]{MortPer}), a compact set $K\subset \SR^n$ is nonpolar if and only if $c_2(K) > 0$.

It is worth mentioning that the ${\bf P} (\tau_K < \infty )$ can take values other than zero or one if the dimension $n$ is at least 3. To see this, consider the ball $\{ \mid x \mid \ < r \}$. Then the probability to hit that ball starting from $x \not \in \{ \mid x \mid \ < r \} $ is given by $ (\frac{r}{\mid x \mid }) ^ {n-2} $ (see e.g. \cite[p.56]{PS}), where $n \geq 3$ is the dimension of the space.

\medskip 

\subsection{Condenser capacity}\label{CC}
We now discuss a related topic, the capacity of a condenser. A {\it condenser} is a pair $(D,K)$, where $D$ is a region and $K \subset D$ is a compact set, both in $\RR^n$. The {\it capacity} of $(D,K)$ is defined to be the infimum of the Dirichlet integral
$$
{\rm cap}(D,K) = \int_{G \backslash D} |\nabla u|^2 m_n(dx),
$$
where the infimum is taken over all smooth functions $u$ with $u=1$ on $K$ and $u=0$ on $\partial D$. This quantity has a number of nice properties that will be relevant to us. For instance, Dirichlet's principle implies that if a minimizer exists then it must be harmonic. When the minimizer does exist it can be interpreted in terms of Brownian motion by $u(z) = {\bf P}^z(B_{T_{D \backslash K}} \in K)$. Furthermore, this type of capacity in two dimensions can be shown to be conformally invariant, in the sense that if $f$ is a conformal map from $D$ onto another domain $D'$, then ${\rm cap}(D,K) = {\rm cap}(D', f(K))$. 
The condenser capacity ${\rm cap}(D,K)$ is also known as the {\it Green capacity} of $K$ with respect to $D$.
See \cite{AV}, \cite{zori}, \cite{PS}  for more on this topic.

\medskip

\subsection{Equilibrium measure}\label{EM}
Let $K$ be a compact subset of a domain $D$ in $\SR^n$. We assume that $D$ possesses a finite Green's function $G_D$. Then there exists a unique Borel measure $\mu_{K,D}$ on $K$ such that
\begin{equation}
{\bf P}^x(\tau_K<T_D)=\int_K G_D(x,y)\,\mu_{K,D}(dy), \;\;\;x\in D.
\end{equation}
This is the {\it equilibrium measure} of $K$ with respect to $D$. Its total measure is equal to the capacity of the condenser $(D,K)$:
$\mu_{K,D}(K)={\rm cap}(D,K)$; see \cite[Chapter 6]{PS}.

\medskip

\subsection{Fundamental frequency}\label{FF}
If $D$ is a planar domain, its {\it fundamental frequency} is given by
\begin{equation}\label{ff}
\lambda(D)=\inf_\phi\frac{\int_D|\nabla \phi|^2\,dm_n}{\int_D\phi^2\,dm_n},
\end{equation}
where the infimum is taken over all smooth functions $\phi$ with compact support in $D$. If the Laplacian has a sequence of Dirichlet eigenvalues on $D$ (e.g. when $m_n(D)<\infty$), then $\lambda(D)$ represents the first eigenvalue. We note, however, that $\lambda(D)$ is defined by (\ref{ff}) even when there are no eigenvalues, and that $\lambda(D)$ may be equal to zero. 

\medskip

The next theorem (see e.g. \cite[\S 3.1]{szi}) gives the connection between the principal Dirichlet eigenvalue and Brownian motion which we will exploit. 

\medskip
  
{\bf Theorem C.}
{\it If $D$ is a domain in $\mathbb C$, then  for every $z\in D$,
\begin{equation}
\lim_{t\to +\infty}\frac{2}{t}\log\frac{1}{{\bf P}^z(T_D>t)}=\lambda(D).
\end{equation} 	
}


The following theorem provides a characterization of the disk as an extremal representative of the class of Schlicht domains; recall that these are images of the unit disk under conformal maps $f$ with $f(0)=0$ and $f^\prime(0)=1$. For a proof, see \cite[\S 5.8]{PoSz} and \cite{He}.
   
   \medskip 
   
{\bf Theorem D.}
   	{\it If $D$ is a Schlicht domain then
   	$\lambda(D)\leq \lambda(\mathbb D)$, with equality if and only if $D=\mathbb D$.
   }

\medskip

Armed with this large assortment of tools, we can now tackle the proofs of our theorems.


\section{{\bf Proofs of Theorems \ref{T1} and \ref{T3}}}  \label{pfT1}

The proof of the theorem differs in details depending on whether the dimension $n$ satisfies $n=2$ or $n \geq 3$. We will give a complete proof in the case $n=2$, and then indicate how the proof must be adjusted when $n \geq 3$. The heart of the proof is contained in the following lemma. 

\begin{lemma}\label{L1} Let $K$ be a compact set in the plane with $0\notin K$ and ${\rm c}_2(K)>0$. Let $a$ be a regular point of $K$.  Let $\delta>0$. There exist  positive constants $C=C(K,\delta)$ and $T=T(K,\delta)$ such that for every $t\in (0,T)$,
\begin{equation}\label{L1e1}
{\bf P}^0(\tau_K<t)\geq C\,e^{-(|a|+\delta)^2/(2t)}.
\end{equation}
\end{lemma}
\proof
We will use a modification of trick taken from \cite[Proof of Lemma 3.6]{GS}. 

Set $L:=K\cap\overline{D(a,\delta/5)}$ and 
 $\Omega:=D(0,3|a|+\delta)$. We will eventually stop the Brownian motion upon exiting $\Omega$, which allows to use the equilibrium measure. Figure 2 may help the reader understand this setup.

  
\begin{figure}
\begin{centering}
\includegraphics[width=10cm,height=10cm,keepaspectratio]{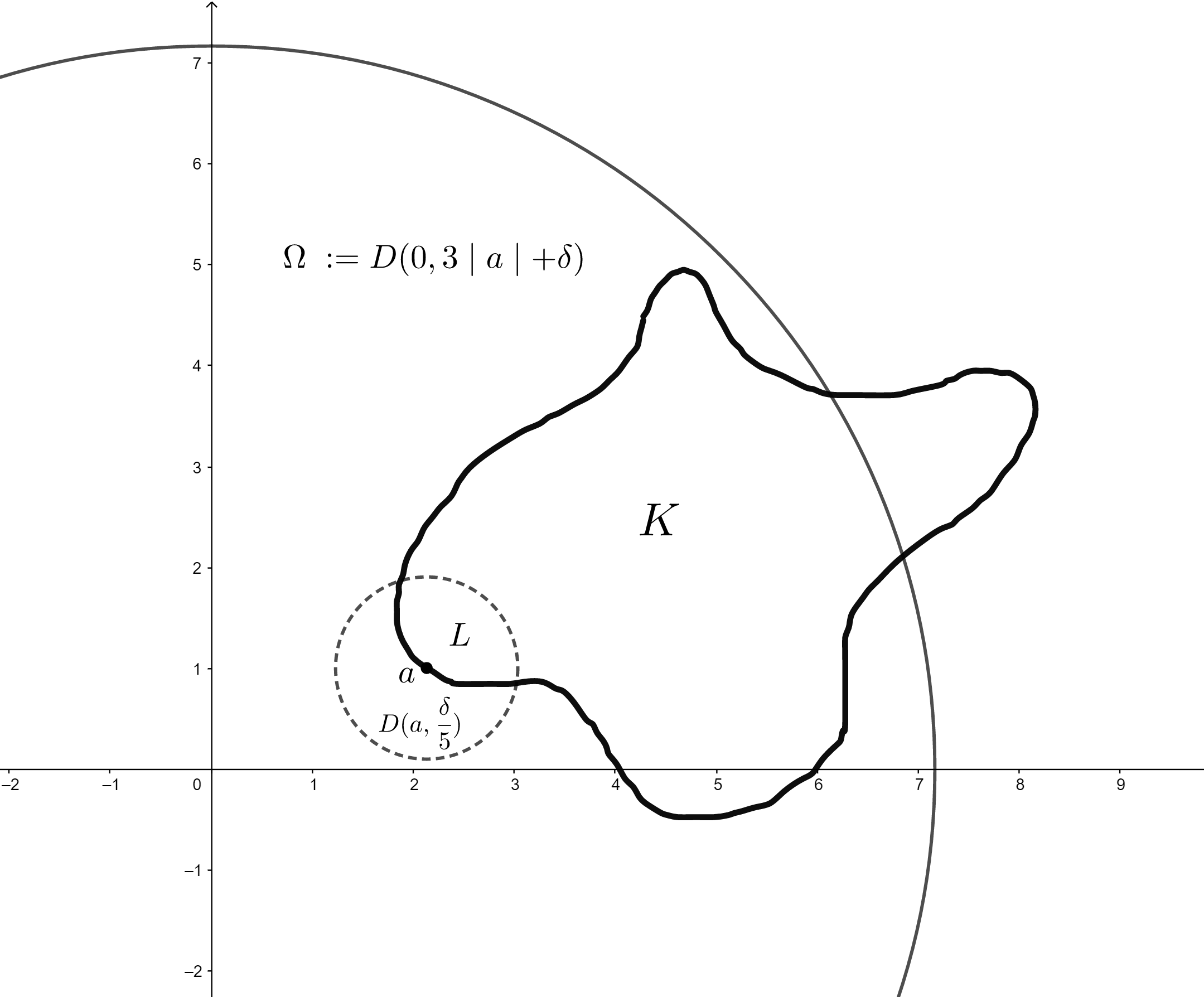}
\par\end{centering}
\caption{$K, L,$ and $\Omega$.}
\end{figure} 
 
 Note that $a$ is a regular point of $L$ and  $c_2(L)>0$. Since $L\subset K$, for every $t>0$ we have
 \begin{eqnarray}\label{L1p4}
{\bf P}^0(\tau_L<t)\leq {\bf P}^0(\tau_K<t) 
\end{eqnarray}

  Observe also that for $t>0$,
\begin{eqnarray}\label{L1p5}
{\bf P}^0(\tau_L<t)&\geq & {\bf P}^0(\tau_L<t\;\;\;\hbox{and}\;\;\;\tau_L<T_\Omega) \\
&=& {\bf P}^0(\tau_L<T_\Omega)- {\bf P}^0(\tau_L\geq t\;\;\;\hbox{and}\;\;\;\tau_L<T_\Omega). \nonumber
\end{eqnarray}

By \cite[Theorem 5.1, p. 190]{PS}, for every $z\in\SC$,
\begin{eqnarray}\label{L1p6}
{\bf P}^z(\tau_L<T_\Omega)&=&\int_L G_\Omega(z,x)\mu_{L,\Omega}(dx) \\
&=&\int_L\int_0^\infty p_\Omega(s,z,x)\, ds\,\mu_{L,\Omega}(dx),\nonumber
\end{eqnarray}
where $G_\Omega$ is the Green function for $\Omega$, $p_\Omega$ is the transition density for Brownian motion killed upon exiting $\Omega$
and $\mu_{L,\Omega}$ is the equilibrium measure of $L$ with respect to $\Omega$.

\medskip

By the (simple) Markov property, the equation (\ref{L1p6}), Fubini's theorem, the domain monotonicity, the semigroup property
of the transition density, and a change of variables, 
\begin{eqnarray}\label{L1p7}
&{}&{\bf P}^0(\tau_L\geq t\;\;\;\hbox{and}\;\;\;\tau_L<T_\Omega) \\
&=&\int_{\Omega\setminus L}p_{\Omega\setminus L}(t,0,z)\;{\bf P}^z(\tau_L<T_\Omega)\;m_n(dz) \nonumber \\
&=& \int_{\Omega\setminus L}p_{\Omega\setminus L}(t,0,z)\int_L\int_0^\infty p_\Omega(s,z,x)ds\;\mu_{L,\Omega}(dx)\;m_n(dz)  \nonumber \\
&=& \int_L\int_0^\infty \int_{\Omega\setminus L} p_{\Omega\setminus L}(t,0,z)  p_\Omega(s,z,x)m_n(dz)\;ds\;\mu_{L,\Omega}(dx) \nonumber \\
&\leq & \int_L\int_0^\infty \int_{\Omega} p_{\Omega}(t,0,z)  p_\Omega(s,z,x)m_n(dz)\;ds\;\mu_{L,\Omega}(dx) \nonumber \\
&=& \int_L\int_0^\infty  p_{\Omega}(t+s,0,x) \,ds\;\mu_{L,\Omega}(dx) \nonumber \\
&=& \int_L\int_t^\infty  p_{\Omega}(s,0,x)ds\;\mu_{L,\Omega}(dx). \nonumber
\end{eqnarray}

Combining (\ref{L1p4}), (\ref{L1p5}),  (\ref{L1p6}), (\ref{L1p7}), we obtain
\begin{eqnarray}\label{L1p8}
&{}& {\bf P}^0(\tau_K<t) \\
&\geq& \int_L\int_0^\infty  p_{\Omega}(s,0,x)ds\;\mu_{L,\Omega}(dx)- \int_L\int_t^\infty  p_{\Omega}(s,0,x)ds\;\mu_{L,\Omega}(dx) \nonumber \\
&=&  \int_L\int_0^t  p_{\Omega}(s,0,x)ds\;\mu_{L,\Omega}(dx).\nonumber
\end{eqnarray}
We will use the fact that for the disk $\Omega$, the transition density $p_\Omega(s,0,x)$ is a decreasing function of $|x|$ (Lemma \ref{L0}).
We will  also use the fact that $\mu_{L,\Omega}(L)$ is equal to the Green capacity (or condenser capacity) ${\rm cap}(\Omega,L)$, and obtain:
\begin{eqnarray}\label{L1p9}
\int_L\int_0^t  p_{\Omega}(s,0,x)ds\;\mu_{L,\Omega}(dx) 
&\geq &  \int_L\int_0^t  p_{\Omega}(s,0,|a|+\delta/5)ds\;\mu_{L,\Omega}(dx)  \nonumber \\
&=& {\rm cap}(\Omega,K)\int_0^t p_\Omega(s,0,|a|+\delta/5)\,ds. 
\end{eqnarray}

\medskip

Now we use ``the principle of not feeling the boundary" (see Subsection \ref{KBM}) and find a positive $T_1=T_1(K,\delta)$ such that
\begin{equation} \label{L1p10}
p_\Omega(s,0,|a|+\delta/5)\geq \frac{1}{2}p(s,0,|a|+\delta/5),\;\;\;\;\;s\in (0,T_1).
\end{equation}	
By (\ref{L1p8}), (\ref{L1p9}), and (\ref{L1p10}), for $t\in (0,T_1)$,
\begin{eqnarray}\label{L1p11}
{\bf P}^0(\tau_K<t) 
&\geq & \frac{{\rm cap}(\Omega,L)}{2}\;\int_0^t p(s,0,|a|+\delta/5)ds  \\
&=&  \frac{{\rm cap}(\Omega,L)}{2}\;\int_0^t \frac{1}{2\pi s}e^{-(|a|+\delta/5)^2/(2s)}\;ds. \nonumber
\end{eqnarray}

By elementary calculus, there exists a positive number $T<T_1$ such that for every $t\in (0,T)$,
\begin{equation} \label{L1p12}
\int_0^t \frac{1}{2\pi s}e^{-(|a|+\delta/5)^2/(2s)}\;ds\geq \frac{1}{2}e^{-(|a|+\delta)^2/(2t)}.
\end{equation}	
Also, by \cite[Lemma 1]{B}, and denoting by $L^*$ the disk $D(0,c_2(L))$, we get
\begin{equation} \label{L1p13}
{\rm cap}(\Omega,L)\geq {\rm cap}(\Omega,\overline{L^*})=\left ( \log\frac{3|a|+\delta}{{\rm c}_2(L)}\right )^{-1},
\end{equation}	
which is a quantity that depends only on $K$ and $\delta$.

\medskip

We combine (\ref{L1p11}), (\ref{L1p12}), (\ref{L1p13}) to conclude that for $t\in (0,T)$,
\begin{eqnarray}\label{L1p14}
{\bf P}^0(\tau_K<t) 
\geq  C\;e^{-(|a|+\delta)^2/(2t)}.
\end{eqnarray}	
\qed

\bigskip

{\it Proof of Theorem \ref{T1}}

We start with some reductions. First, we assume without loss of generality that 
${\rm d}((\partial W)^r)=1$. Then for every $t>0$,
\begin{equation} \label{T1p1}
{\bf P}^0(T_W<t)\leq {\bf P}^0(T_\SD<t),
\end{equation}	
and so we may assume that $W=\SD$. 

Let $a\in (\partial U)^r$ be a point with $|a|<1$. Choose a disk $D(a,r)$ with $0<r<1-|a|$, and
set $K:=(\SC\setminus U)\cap \overline{D(a,r)}$. Note that $K$ has positive logarithmic capacity and that 
$U\subset \SC\setminus K$. Therefore,
\begin{equation} \label{T1p2}
{\bf P}^0(T_U<t)\geq {\bf P}^0(T_{\SC\setminus K}<t)={\bf P}^0(\tau_K<t).
\end{equation}	
So we may assume that $U=\SC\setminus K$, where $K$  is a compact subset of $\SD$ with ${\rm c}_2(K)>0$.

\medskip

By Lemma \ref{L1}, for every $\delta>0$, there exist positive constants $C,T$ depending on $K$ and $\delta$
such that for every $t\in (0,T)$,
\begin{equation} \label{T1p3}
{\bf P}^0(\tau_K<t)\geq  C\;e^{-(|a|+\delta)^2/(2t)},
\end{equation}	
where $a$ is a regular point of $K$. By taking $\delta$ small enough, we may assume that $|a|+\delta<1$. On the other hand, an ingenious argument due to McConnell \cite{MAC} shows that, for all $t>0$ and
all positive integers $m\geq 3$,
\begin{equation} \label{macest}
{\bf P}(T_\SD<t)\leq c(m)\;e^{-\frac{\cos^2(\pi/m)}{2t}}.
\end{equation}
By fixing $m$ large enough, we see that for any $\epsilon >0$ we can find a constant $C_\epsilon > 0$ such that for all $t>0$ we have
\begin{equation*}
{\bf P}(T_\SD\leq t)<C_\epsilon\;e^{-\frac{(1-\epsilon)}{2t}}.
\end{equation*}
Thus Lemma \ref{L1} gives
\begin{eqnarray}
\lim_{t\to 0+}\frac{{\bf P}(T_U<t)}{{\bf P}(T_W<t)} &\geq & \lim_{t\to 0+}\frac{{\bf P}(\tau_K<t)}{{\bf P}(T_\SD<t)} \\
&\geq & \lim_{t\to 0+}\frac{C(K,\delta)\;\exp \left [-\frac{(|a|+\delta)^2}{2t}\right  ]}{C_\epsilon\;\exp\left [-\frac{1-\epsilon}{2t}\right ]}. \nonumber
\end{eqnarray}
Choose $\epsilon$ and $\delta$ small enough so that $(|a| +\delta)^2<1-\epsilon$. Then the limit above is equal to $\infty$ and this completes the proof of the result in two dimensions.
 
\medskip

The proof in three or higher dimensions follows along exactly the same lines, with several necessary modifications, which we now indicate. In $n \geq 3$ dimensions, we will use the Newtonian capacity rather than the logarithmic. All other concepts used carry over directly to the higher dimensions with no real change, except that some of the estimates have to be changed. In particular, inequality (\ref{L1p13}) is specific to two dimensions; it is sufficient to replace it with \cite[Cor. 1]{zori}. Furthermore, McConnell's estimate (\ref{macest}) is also specific to two dimensions, but it can  be replaced by \cite[Cor. 3.4]{sera}. The result follows then as before.
\qed


\bigskip         \bigskip

{\it Proof of Theorem \ref{T3}}

 Suppose that $D$ is a Schlicht domain other than $\mathbb D$. By Theorems C and D (in Subsection \ref{FF}),
   \begin{eqnarray}
   \lim_{t\to+\infty}\left (\frac{{\bf P}^0(T_D>t)}{{\bf P}^0(T_{\mathbb D}>t)}\right )^{2/t} &=&
   \lim_{t\to+\infty}\exp \left [\frac{2}{t}\log\frac{{\bf P}^0(T_D>t)}{{\bf P}^0(T_{\mathbb D}>t)} \right ] \nonumber\\
   &=& \exp(\lambda(\mathbb D)-\lambda(D))>1.
   \end{eqnarray}
   It follows that there exists $t_o>0$ such that for every $t>t_o$, 
   $$
   {\bf P}^0(T_{D}>t)>{\bf P}^0(T_{\mathbb D}>t),
   $$ 
   which is equivalent to (\ref{T1e1new}).
 \qed

\bibliographystyle{amsplain}

\end{document}